\documentclass[a4paper, 11pt]{amsart}

\usepackage{amsthm}
\usepackage{amsmath}
\usepackage{amsfonts}
\usepackage{amssymb}

\theoremstyle{plain}
\newtheorem{theorem}{\textbf{Theorem}}[section]

\newtheorem{corollary}[theorem]{\textbf{Corollary}}
\newtheorem{proposition}[theorem]{\textbf{Proposition}}

\newtheorem*{Freiman-3k-4}{\textbf{Freiman $3k-4$ Theorem}}
\newtheorem*{Freiman-3k-3}{\textbf{Freiman $3k-3$ Theorem}}

\newcommand{\Z}{\mathbb{Z}}

\newcommand{\be}{\begin{equation}}
\newcommand{\ee}{\end{equation}}

\newcommand{\diam}{\mathsf{diam\,}}

\newcommand{\ber}{\begin{eqnarray}}
\newcommand{\eer}{\end{eqnarray}}
\newcommand{\nn}{\nonumber}
\newcommand{\und}{\;\;\;\mbox{ and }\;\;\;}

\newcommand{\gcds}{\gcd^*}

\begin{document}

\title{Long Arithmetic Progressions in Sets with Small Sumset}

\author{Itziar Bardaji}
\author{David J. Grynkiewicz$^*$}

\thanks{$*$ supported by FWF project number M1014-N13}

\address{Institut de Rob\`{o}tica i Inform\`{a}tica Industrial, Universitat Polit\`{e}cnica de Catalunya, 08034-Barcelona, Espanya}\email{ibardaji@iri.upc.edu}

\address{Institut f\"ur Mathematik und Wissenschaftliches Rechnen,
Karl-Franzens-Universit\"at Graz,
Heinrichstra\ss e 36,
8010 Graz, Austria} \email{diambri@hotmail.com}

\begin{abstract}
 Let $A,\,B\subseteq \Z$ be finite, nonempty subsets with $\min A=\min B=0$, and let $$\delta(A,B)=\left\{
  \begin{array}{ll}
    1 & \hbox{if } A\subseteq B, \\
        0 & \hbox{otherwise.}
  \end{array}
\right.$$  If $\max B\leq \max A\leq |A|+|B|-3$ and
\be\label{card-hypothII}|A+B|\leq |A|+2|B|-3-\delta(A,B),\ee  then we show $A+B$ contains an arithmetic progression with
difference $1$ and length $|A|+|B|-1$.

As a corollary, if \eqref{card-hypothII} holds, $\max(B)\leq \max(A)$ and either $\gcd(A)=1$ or else  $\gcd(A+B)=1$ and $|A+B|\leq 2|A|+|B|-3$, then $A+B$ contains an arithmetic progression with
difference $1$ and length $|A|+|B|-1$.
\end{abstract}

\maketitle

\section{Introduction}

For a subset $A\subseteq \Z$, we let $\diam A=\max A-\min A$ denote
its diameter and $|A|$ its cardinality. We let $\gcds(A)=\gcd(A-a_0)$, where $a_0\in A$ and $\gcd$ denotes the greatest common divisor. For $A,\,B\subseteq\Z$, their sumset is the set of all sums of
one element from $A$ and another from $B$: $$A+B=\{a+b\ :\ a\in A,\,
b\in B\}.$$ Also, define $$\delta(A,B)=\left\{
  \begin{array}{ll}
    1 & \hbox{if } x+A\subseteq B, \hbox{ for some } x\in \Z\mbox{,} \\
        0 & \hbox{otherwise.}
  \end{array}
\right.$$

The study of the structure of subsets with small sumset has a rich tradition (see \cite{natbook} and \cite{taobook} for two texts on the subject). One classical result is the $(3k-4)$-Theorem of Freiman \cite{Freiman-3k-4-sym} \cite{Freiman-Monograph} \cite{natbook} \cite{taobook}, which  states that if a set $A$ of integers satisfies $\gcds(A)=1$ and
\be\label{3k-4bound}|A+A|\le 3|A|-4,\ee then the diameter of $A$ is at most $2|A|-4$. In other words, $A$ is an interval with at most $|A|-3$ holes. Various generalizations to distinct summands were later found \cite{Freiman-3k-4-distinct} \cite{Lev-Smel-3k-4} \cite{Stanchescu-3k-4} \cite{3k-2-distinct-summands}. The latest result is from \cite{3k-2-distinct-summands}, and shows that if $\diam (A)\geq \diam(B)$, $\gcds(A)=1$ and \be\label{3k-4bound-distinct} |A+B|=|A|+|B|-1+r\leq |A|+2|B|-3-\delta(A,B),\ee then $\diam(A)\leq |A|+r-1$ and $\diam(B)\leq \min\{|A|,\,|B|\}+r-1$.

However, as has later become apparent, knowing that there are only a small number of holes is not always sufficient. In part, this is because there are many subsets of small diameter that nonetheless have large sumset. Working through examples, one quickly finds that, informally speaking, it is much more difficult for the holes in a subset $A$ with small sumset (and correspondingly the holes in $A+A$ as well) to occur in the interior of the set than near the boundary (namely, near the maximum or minimum element). However, there have been few results satisfyingly embodying this idea.

One such result occurred as a lemma in a pair of papers of J. Deshouillers and V. Lev characterizing large sum-free set over $\Z/p\Z$ \cite{lev-apI} \cite{des-lev}. It's main consequence stated that if $\diam(A)<\frac{3}{2}|A|-1$, then $A-A$ contains an interval of length $2|A|-1$ \cite[Lemma 3]{lev-apI}. (Some other similar results were also given.) Those familiar with Freiman's Theorem (see \cite{natbook} or \cite{taobook}) may also recall that the existence of large multi-dimensional progressions plays an important role in its proof. Very recently, G. Freiman gave a very precise estimate for the length of an arithmetic progression that can be found in $A+A$ when the sumset is so small as to satisfy \eqref{3k-4bound}, showing that there is always one of length at least $2|A|-1$. The example $$A=\{0,1,2,\ldots,k-r-1,k-r+1,k-r+3,\ldots,k-r-1+2r\},$$ for $r=0,1,\ldots,k-3$, shows the bound on the arithmetic progression length to be best possible, while the example $$A=\left\{1,2,\ldots,\left\lceil\frac{k}{2}\right\rceil\right\}\cup\left\{x+1,x+2,\ldots,x+\left\lfloor \frac{k}{2}\right\rfloor\right\},$$ for  $x\geq k+1$, shows that the assumption $|A+A|\leq 3|A|-4$ from \eqref{3k-4bound} is needed. The paper of Freiman also delved into the issue of where the holes could occur in $A$, but the other structural information is derivable from the bound on the length of the arithmetic progression.

The goal of our paper is to extend this result of Freiman to pairs of distinct summands $A$ and $B$.

\begin{theorem}\label{mainresult-intro-formulation} Let $A,\,B\subseteq \Z$ be nonempty subsets with $\diam B\leq \diam A\leq |A|+|B|-3$ and
\be\label{card-hypoth}|A+B|\leq |A|+2|B|-3-\delta(A,B).\ee  Then $A+B$ contains an arithmetic progression with
difference $1$ and length $|A|+|B|-1$.\end{theorem}

Using the previously mentioned structural result for small sumsets from \cite{3k-2-distinct-summands}, we obtain the following immediate corollary. Note $\gcds(A)>1$ and $\gcds(A+B)=1$ trivially implies $|A+B|\geq 2|A|+|B|-2$; see \cite[Theorem 1.2]{natbook}.

\begin{corollary}\label{main-cor}Let $A,\,B\subseteq \Z$ be nonempty subsets with $\diam B\leq \diam A$ and $\gcds(A+B)=1$. If  \be\label{tee2}|A+B|\leq |A|+2|B|-3-\delta(A,B)\ee and either $\gcds(A)=1$ or  \be\label{tee}|A+B|\leq 2|A|+|B|-3,\ee then
$A+B$ contains an arithmetic progression with
difference $1$ and length $|A|+|B|-1$.
\end{corollary}

We remark that the condition $\gcds(A+B)=1$ is simply a normalization hypothesis; if instead $\gcds(A+B)=d$ and \eqref{tee2} and \eqref{tee} hold, then the difference of the arithmetic progression becomes $d$.

During the course of the proof of Theorem \ref{mainresult-intro-formulation}, the structural consequences concerning the location of holes and such will become apparent in the series of propositions and definitions leading up to the proof of Theorem \ref{mainresult-intro-formulation}. The paper concludes with a few additional remarks.

\section{Long Arithmetic Progressions}

Throughout this section, we assume $A,\,B\subseteq \Z$ are finite,
nonempty subsets normalized so that
\be\label{minABis0} \min A=\min B=0,\ee
and
 with \ber \label{MN-def} M=\max
A\;\;\;&\mbox{ and }&\;\;\; N=\max B,\\\label{M-bigger-N} M&\geq&
N,\\\label{card-assumption} |A+B|&=&|A|+|B|-1+r,\eer so that $A$ is
assumed to be the set with larger (or equal) diameter. As the problem is translation invariant, there is no loss of generality when assuming \eqref{minABis0}. Note, in view of \eqref{minABis0} and \eqref{M-bigger-N}, that \be\label{remarkk}\delta(A,B)=1\mbox{ if and only if }A\subseteq B.\ee

For $a,\,b\in \Z$, we define $[a,b]:=\{x\in \Z\mid a\leq x\leq b\}\subseteq \Z$. For a set $X$ and an interval $[a,b]\subseteq \Z$, the number of holes of $X$ in $[a,b]$ is denoted by
$$h_{X,[a,b]}=|[a,b]\setminus X|.$$ When $[a,b]$ is the default interval $[\min X, \max X]$, we skip reference to the
interval, that is, $$h_X=h_{X,[\min X,\max X]},$$ and when we refer to a hole in $X$ without reference to an interval, we simply mean an element $x\in [\min X,\max X]\setminus X$.

Observe, in view of \eqref{minABis0}, \eqref{MN-def} and \eqref{card-assumption}, that \ber
\label{expression-for-M} M&=&|A|+h_A-1,\\N&=&|B|+h_B-1,\label{expression-for-N}\\
\label{expression-for-h-AB}h_{A+B}&=&M+N+1-|A+B|=h_A+h_B-r.\eer
Also remark that, using \eqref{expression-for-M}, we can rewrite the condition $\diam A\leq |A|+|B|-3$ in Theorem \ref{mainresult-intro-formulation} as $h_A\leq |B|-2$ and the condition \eqref{card-hypoth} as $r\le |B|-2-\delta(A,B)$.

\bigskip

\begin{proposition}\label{prop-neutral-zone} If $h_A\leq |B|-1$,
then \be\label{neutral-zone-eq} [N,M]\subseteq A+B.\ee
\end{proposition}

\begin{proof}
 Let $x\in [N,M]$. Thus
 $$(x,0),\,(x-1,1),\,\ldots,(x-N,N)$$ are all representations $(a,b)$ of
 $x=a+b$ with $a\in [0,M]$ and $b\in [0,N]$. If $x\notin A+B$, then
each of these $N+1$ pairs must either have the first element missing
from $A$ or the second element missing from $B$, whence $h_A+h_B\geq
N+1=|B|+h_B$ (in view of \eqref{expression-for-N}). But this
contradicts $h_A\leq |B|-1$.
\end{proof}

\bigskip

In view of Proposition \ref{prop-neutral-zone}, we see that all
holes in $A+B$ lie in one of the disjoint intervals $[0,N-1]$ or
$[M+1,M+N]$. We refer to them as \emph{left} and \emph{right} holes, respectively.


 Since $0+B\subseteq A+B$ and $0+A\subseteq A+B$, if $x\in [0,N-1]$ is a hole
in $A+B$, then $x$ must also be a hole in $B$ and in $A$. Likewise, since $M+B\subseteq A+B$, if $x+M\in [M+1,M+N]$ is a hole
in $A+B$, then $x$ must also be a hole in $B$, and since
$N+A\subseteq A+B$, if $x+N\in [M+1,M+N]$ is a hole in $A+B$, then
$x$ must also be a hole in $A$.

In view of these observations, we make some definitions.

\begin{itemize}
 \item Holes $x\in
[0,N-1]\setminus B$ for which $x\notin A+B$ remains a hole in $A+B$
 are called
\emph{left stable} holes in $B$.
\item Holes $x\in [0,N-1]\setminus A$ for which $x\notin A+B$ are called \emph{left stable} holes in $A$.
\item Holes $x\in
[1,N]\setminus B$ for which $x+M\notin A+B$ are called \emph{right stable} holes in $B$.
\item Holes $x\in [M-N+1,M]\setminus A$
for which $x+N\notin A+B$ are called
\emph{right stable} holes in $A$.
\end{itemize}

A \emph{stable} hole in $A$ is
one which is either right or left stable, and likewise for $B$. All
other holes (in either $A$ or $B$) are called \emph{unstable}. We
let $h_A^s$ and $h_B^s$ denote the respective number of stable
holes in $A$ and $B$, and we let $h_A^u$ and $h_B^u$ denote the
respective number of unstable holes in $A$ and $B$.

This classification of holes into ones which contribute to a hole
present in $A+B$ (the stable ones) and those which do not contribute
to any hole in $A+B$ (the unstable ones) will prove to be a very
useful perspective.
Note that a pair of stable holes
$x_A$ and $x_B$, one from $A$ and one from $B$, can be associated to each hole $x\notin A+B$: Indeed, if $x\in
[0,N-1]$ is a left hole in $A+B$, then $x$ must come from a left stable hole both in $B$ and
$A$, i.e., $$x_B=x\notin B\und x_A=x\notin A$$ are both left stable
holes. On the other hand, if $x\in [M+1,M+N]$ is a right hole in $A+B$, then it must come
from a right stable hole both in $B$ and $A$, i.e., $$x_B=x-M\notin
B\und x_A=x-N\notin A$$ are both right stable holes.

We will later see that these mappings are invertible, i.e., that
$x_A=y_A$ for holes $x,\,y\notin A+B$ implies $x=y$, and likewise
$x_B=y_B$ implies $x=y$. However, next we prove a very important
proposition---the key observation used in the proof---which shows
that if we have a left hole $x\notin A+B$,
then there must be many holes in $A\cap [0,x]$ and $B\cap [0,x]$,
with an analogous statement holding for right holes.

\bigskip

\begin{proposition}\label{Key-lemma} If $x\in [0,N]\setminus (A+B)$,
then \be\label{left-hole-bound} h_{A,[0,x]}+h_{B,[0,x]}\geq x+1.\ee
If $x+M\in [M,M+N]\setminus (A+B)$, then \be\label{right-hole-bound}
h_{A,[x+M-N,M]}+h_{B,[x,N]}\geq N-x+1.\ee
\end{proposition}

\begin{proof}
The proof is analogous to that of the previous proposition.  If
$x\in [0,N]$, then $$(x,0),\,(x-1,1),\,\ldots,(0,x)$$ are all
representations $(a,b)$ of
 $x=a+b$ with $a\in [0,M]$ and $b\in [0,N]$ (in view of \eqref{M-bigger-N}). If $x\notin A+B$, then
each of these $x+1$ pairs must either have the first element missing
from $A$ or the second element missing from $B$, whence
\eqref{left-hole-bound} follows. The argument for when $x+M\in
[M,M+N]$ is analogous, considering instead
$$(M,x),\,(M-1,x+1),\,\ldots,(x+M-N,N).$$
\end{proof}

\bigskip

Next, we show that no hole in $B$ can be both left and right stable.

\bigskip

\begin{proposition}\label{lem-left-right-B} Let $x\in [1,N]\setminus
B$. If $h_A\leq |B|-2$, then 
either $x\in A+B$ or $x+M\in A+B$.\end{proposition}

\begin{proof}
If both $x\notin A+B$ and $x+M\notin A+B$, then applying both cases
of Proposition \ref{Key-lemma} yields \ber\nn
N+2=(x+1)+(N-x+1)&\leq&
h_{A,[0,x]}+h_{B,[0,x]}+h_{A,[x+M-N,M]}+h_{B,[x,N]}\\\nn&\leq&
h_A+h_B+2,\eer where the second inequality follows by
\eqref{M-bigger-N}. Now applying \eqref{expression-for-N} yields
$h_A\geq |B|-1$, contrary to assumption.
\end{proof}

\bigskip

The following shows there are also no holes in $A$ which are both
left and right stable.

\bigskip

\begin{proposition}\label{lem-left-right-A} Let $x\in [0,M]\setminus
A$. If $h_A\leq |B|-2$, then 
either $x\in A+B$ or $x+N\in A+B$.\end{proposition}

\begin{proof}
If both $x\notin A+B$ and $x+N\notin A+B$, then Proposition \ref{prop-neutral-zone} implies $x\in [M-N+1,N-1]$, whence applying both cases
of Proposition \ref{Key-lemma} yields \ber\nn
M+2=x+1+(N-(x-M+N)+1)&\leq&
h_{A,[0,x]}+h_{B,[0,x]}+h_{A,[x,M]}+h_{B,[x-(M-N),N]}\\\nn &\leq&
h_A+1+h_B+M-N+1, \eer where the second inequality follows in view of
\eqref{M-bigger-N}. Now applying \eqref{expression-for-N} yields
$h_A\geq |B|-1$, contrary to assumption.
\end{proof}

\bigskip

In view of Propositions \ref{prop-neutral-zone},
\ref{lem-left-right-B} and \ref{lem-left-right-A} (and the relevant
stability definitions), we see that, when $h_A\leq |B|-2$,
\ber\label{stab-eqs} h_A^s=h_B^s=h_{A+B}&=&h_A+h_B-r\\\label{unstab-eq-A} h_A^u&=&r-h_B\\
h_B^u&=&r-h_A\label{unstab-eq-B},\eer where \eqref{stab-eqs} uses
\eqref{expression-for-h-AB}, and where \eqref{unstab-eq-A} and
\eqref{unstab-eq-B} follow from \eqref{stab-eqs} by using the
identities $h_B=h_B^u+h_B^s$ and $h_A=h_A^u+h_A^s$; moreover,
$x_A=y_A$ or $x_B=y_B$ implies $x=y$ for holes $x,\,y\in
[0,M+N]\setminus (A+B)$, as previously alluded. Note $x_B=x_A$ when
$x$ is a left hole, and that $x_A=x_B+(M-N)$ when $x$ is a right
hole.

The next proposition is the trickiest part of the proof, showing
that all left stable holes precede all right stable holes, so there
is no overlap.

\bigskip

\begin{proposition}\label{No-Overlap} Suppose $h_A\leq |B|-2$ and $r\leq |B|-2-\delta(A,B)$.
If $x_B\in [0,N]\setminus B$ is a left stable hole and $y_B\in
[0,N]\setminus B$ is a right stable hole, then $x_B<y_B$. Likewise,
if $x_A\in [0,M]\setminus A$ is a left stable hole and $y_A\in
[0,M]\setminus A$ is a right stable hole, then $x_A<y_A$.
\end{proposition}

\begin{proof} If $x_A\in [0,M]\setminus A$ is a left stable hole,
$y_A\in [0,M]\setminus A$ is a right stable hole and $x_A\geq y_A$,
then $x_B=x_A\in [0,N]\setminus B$ is a left stable hole and
$y_B=y_A-(M-N)\in [0,N]\setminus B$ is a right stable hole with
$x_B\geq y_B$, in view of $x_A\geq y_A$ and \eqref{M-bigger-N}.
Therefore we see that it suffices to prove the first assertion in
the proposition, as the second is an immediate consequence.

To that end, assume $x_B\in [0,N]\setminus B$ is a left stable hole
and $y_B\in [0,N]\setminus B$ is a right stable hole with $x_B> y_B$. Note that $x_B= y_B$ cannot hold in view of Proposition
\ref{lem-left-right-B}. Moreover, assume $x_B$ and $y_B$ are chosen
minimally, meaning that there are no stable holes $z\in
[y_B+1,x_B-1]\setminus B$.

Applying both cases of Proposition \ref{Key-lemma} to $x_B$ and
$y_B+M$, respectively, we find that
\ber\nn |B|+h_B +(x_B-y_B+1)&=&(x_B+1)+(N-y_B+1)\\ \nn &\leq&
h_{A,[0,x_B]}+h_{B,[0,x_B]} +h_{A,[y_B+M-N,M]}+h_{B,[y_B,N]} \\
&\leq & h_A+h_B+h_{A,[y_B+M-N,x_B]}+h_{B,[y_B,x_B]},\label{antsrun}
\eer where we use \eqref{expression-for-N} for the first equality.

In view of the minimality of $x_B$ and $y_B$,  we see that \be\label{estimate-hB1}h_{B,[y_B,x_B]}\leq
h_B^u+2,\ee
with equality
possible  only if $[y_B+1,x_B-1]$ contains all the
unstable holes in $B$.
We also have the trivial inequality
\be\label{estimate-hB2}h_{B,[y_B,x_B]}\leq x_B-y_B+1.\ee
If $y_B+M-N>x_B$, so
that $h_{A,[y_B+M-N,x_B]}=0$, then \eqref{antsrun} and \eqref{estimate-hB2} imply $h_A\geq
|B|$, contrary to hypothesis. Therefore we may assume $y_B+M-N\leq
x_B$, and now we also have the trivial inequality
\be\label{estimate-hA2}h_{A,[y_B+M-N,x_B]}\leq x_B-y_B+1-(M-N),\ee with equality possible
only if $[y_B+M-N,x_B]$ are all holes in $A$.

Applying the estimates \eqref{estimate-hA2} and
\eqref{estimate-hB1} in \eqref{antsrun} and
using \eqref{expression-for-M}, \eqref{expression-for-N} and
\eqref{unstab-eq-B}, we discover that \be\label{deewt} |A|-2-r\leq
h_B-h_A.\ee

In view of \eqref{M-bigger-N}, \eqref{expression-for-M} and \eqref{expression-for-N}, we have \be\label{onemore}h_B-h_A\leq |A|-|B|,\ee with equality only possible when $M=N$. Combining \eqref{onemore} and \eqref{deewt} yields \be\label{iip}r\geq |B|-2,\ee whence our hypothesis $r\leq |B|-2-\delta(A,B)$ implies that $r=|B|-2$, that $\delta(A,B)=0$, and that equality held in all estimates used to derive \eqref{iip}.

As a result, $\delta(A,B)=0$ and \eqref{remarkk} imply $A\nsubseteq B$; equality in \eqref{onemore} implies $M=N$; and equality in \eqref{deewt} implies equality holds in both \eqref{estimate-hA2} and \eqref{estimate-hB1}, whence $[y_B+M-N,x_B]$ are all holes in $A$ and $[y_B+1,x_B-1]$ contains all the
unstable holes in $B$.

Since $A\nsubseteq B$, it follows that there exists $z\in A$ with $z\notin B$. Since $[y_B+M-N,x_B]$ are all holes in $A$ and $M=N$, it follows that $z\notin [y_B,x_B]$.
Thus, since $[y_B+1,x_B-1]$ contains all the
unstable holes in $B$, it follows that $z\notin B$ is a stable hole in $B$. However, this means that either $z+0\notin A+B$ or $z+M\notin A+B$,
which are both contradictions in view of $z\in A$ and $M=N$, completing the proof.
\end{proof}

\bigskip

We are now ready to finish the proof of Theorem \ref{mainresult-intro-formulation}, which we will follow from the next proposition.

\bigskip

\begin{proposition}\label{long-ap-thm} Suppose $h_A\leq |B|-2$ and $r\leq |B|-2-\delta(A,B)$.
Then $$J:=[e+1,M+c-1]\subseteq A+B,$$ where $e$ is the greatest left
stable hole in $B$ (let $e=-1$ if there are no left stable holes)
and $c$ is the smallest right stable hole in $B$ (let $c=N+1$ if
there are no right stable holes). Moreover, \ber\label{J-is-long}
|J|=M-1+(c-e)&\geq&
|A|+|B|-1+h_{A,[e+1,c+M-N-1]}+h_{B,[e+1,c-1]}\\\nn &\geq&
|A|+|B|-1.\eer
\end{proposition}

\begin{proof}
In view of proposition \ref{No-Overlap}, we have $e<c$. Consequently, by
the definition of stability, and in view of Proposition
\ref{prop-neutral-zone}, we see that $$J:=[e+1,M+c-1]\subseteq
A+B.$$ Note \be\label{elecktra-D}|J|=M-1+(c-e)=|A|+h_A-2+(c-e),\ee
using \eqref{expression-for-M}. It remains to estimate $c-e$.

Applying both cases of Proposition \ref{Key-lemma} to $e$ and $c+M$,
respectively, we find that \be\label{truck}e+1+N-c+1\leq
h_A+h_B-s,\ee where $s=h_{A,[e+1,c+M-N-1]}+h_{B,[e+1,c-1]}$. From \eqref{expression-for-N} and \eqref{truck}, it follows that $$|B|+h_B+1+e-c\leq
h_A+h_B-s,$$ yielding $$c-e\geq |B|+1-h_A+s.$$ Combining the
above estimate for $c-e$ with \eqref{elecktra-D}, we obtain \ber\nn
|J|=|A|+h_A-2+(c-e)&\geq&
|A|+h_A-2+(|B|+1-h_A+s)\\\nn&=&|A|+|B|-1+s\geq |A|+|B|-1,\eer
completing the proof.
\end{proof}

\bigskip

Finally, we complete the proof of Theorem \ref{mainresult-intro-formulation}.

\begin{proof}[Proof of Theorem \ref{mainresult-intro-formulation}] We may w.l.o.g. assume $\min A=\min B=0$. Since $\diam B\leq \diam A\leq |A|+|B|-3$, we have $h_A\leq |B|-2$ in view of \eqref{MN-def}. Since $|A+B|:=|A|+|B|-1+r\leq |A|+2|B|-3-\delta(A,B)$, we have $r\leq |B|-2-\delta(A,B)$. Thus applying Proposition \ref{long-ap-thm} completes the proof.
\end{proof}

\section{Concluding Remarks}

We conclude with some brief remarks, for which we assume the notation of the previous section, particularly concerning Proposition \ref{long-ap-thm}.

First, let us show that all the intermediary work and propositions leading up to Theorem \ref{mainresult-intro-formulation}, save Proposition \ref{Key-lemma}, are easily deduced from Theorem \ref{mainresult-intro-formulation} itself.  If $J=[a,b]\subseteq A+B$ is the arithmetic progression with difference $1$ given by Theorem \ref{mainresult-intro-formulation}, then observe that $$(A\cup [a,b-N])+(B\cup [a,b-M])=A+B.$$  From this observation, the apparently stronger bound given by Proposition \ref{long-ap-thm} is now easily derived from Theorem \ref{mainresult-intro-formulation}.
Additionally, if $J$ does not contain the interval $[N,M]$, then it follows that $$|J|\leq M=|A|+h_A-1\leq |A|+|B|-3-\delta(A,B),$$ contrary to Theorem \ref{mainresult-intro-formulation}. Thus Proposition \ref{prop-neutral-zone} is a consequence of Theorem \ref{mainresult-intro-formulation}. Noting that the hypothesis $h_A\leq |B|-2-\delta(A,B)$ implies (in view of \eqref{M-bigger-N} and \eqref{MN-def}) that $h_B\leq |A|-2-\delta(A,B)$, it is then easily derived from the existence of the long arithmetic progression $J=[a,b]\subseteq A+B$, and a simple calculation, that all stable holes preceding $a$ are left stable and that all stable holes following $b-N$ or $b-M$, respectively for $A$ or $B$, are right stable. So the propositions concerning right and left stable holes also follow from Theorem \ref{mainresult-intro-formulation}. This leaves only Proposition \ref{Key-lemma} as containing information additional to Theorem \ref{mainresult-intro-formulation}, as claimed.

Next, it is important to note that Theorem \ref{mainresult-intro-formulation}/Proposition \ref{long-ap-thm} essentially shows that the sets $A$ and $B$ can be divided into left and right halves with each half behaving independently (with respect to the sumset $A+B$) of the other. Taking the left halves $A_L=A\cap [0,e]$ and $B_L=B\cap [0,e]$ and unioning each with a sufficiently long interval $[e+1,x]$, where $x\geq e+1+r$, results in a pair of subsets whose sumset exhibits the same structural behavior on the left side as for the original sumset $A+B$. (The right side of $A+B$ can be independently studied in a similar manner.)

In general, there are many possibilities for how the holes can be distributed in $A_L$ and $B_L$. However, if one wishes to use holes efficiently, that is, use a large number of holes relative to the maximal bound $r$, then \eqref{unstab-eq-A} and \eqref{unstab-eq-B} show that the number of unstable holes must be small, which helps restrict the possibilities for $A_L$ and $B_L$.

For instance, in the extremal case when there are no unstable holes in either $A$ or $B$, then we must have $A_L=B_L$, and $A_L\cup [e+1,\infty)$ is the complement of the solution set of the Frobenius problem (see \cite{Frob}) for the set $A$ (i.e, $A_L\cup [e+1,\infty)=\bigcup_{h=1}^{\infty} hA$, where $hA=A\underbrace{+\ldots+}_h A$ denotes the $h$-fold sumset). In particular, if $d_1,\,d_2\in A_L$, then the arithmetic progression $\{d_1+id_2\mid i=0,1,2,\ldots,\}$ is contained in $A_L\cup [e+1,\infty)$. In fact, $A_L$ is just the intersection of the multi-dimensional progression $\{i_1d_1+i_2d_2+\ldots+i_ld_l\mid i_j=0,1,2,\ldots\}$ with $[0,e]$, where $A_L=\{0,d_1,d_2,\ldots,d_l\}$.


\begin{thebibliography}{99}

\bibitem{des-lev} J. Deshouillers and V. Lev, A refined bound for sum-free sets in groups of prime order, \emph{Bull. Lond. Math. Soc.},  40  (2008),  no. 5, 863--875.

\bibitem{Freiman-Monograph} G.~A.~Freiman, \emph{Foundations of a structural theory of set
addition}, translated from the Russian, Translations of Mathematical
Monographs, \textbf{37}, American Mathematical Society, Providence,
R. I., 1973.

\bibitem{Freiman-3k-4-sym} G. A. Freiman, The addition of finite sets I, \emph{Izv. Vyss. Ucebn. Zaved.
Matematica} 6 (1959), no. 13, 202--213 (Russian) .

\bibitem{Freiman-3k-4-distinct} G.~A.~Freiman,
Inverse problems of additive number theory VI: On the addition of
finite sets III, \emph{Izv. Vys . S.U.cebn. Zaved.
Matematika} (1962), no. 3 (28), 151--157 (Russian).

\bibitem{Freiman-3k-4w-holes-sym} G. A. Freiman, Inverse Additive Number Theory XI: On the detailed structure of sets with small additive property, preprint.

\bibitem{3k-2-distinct-summands} D. J. Grynkiewicz and O. Serra, The Freiman $3k-2$ Theorem: Distinct Summands, preprint.

\bibitem{Lev-Smel-3k-4} V.~F.~Lev and P.~Y.~Smeliansky,
On addition of two distinct sets of integers, \emph{Acta Arith.}
\textbf{70} (1995), no. 1, 85--91.

\bibitem{lev-apI}
V. Lev, Large sum-free sets in $\Z/p\Z$,  \emph{Israel J. Math.},  154  (2006), 221--233.

\bibitem{natbook} M.~Nathanson, \emph{Additive Number Theory:
Inverse Problems and the Geometry of Sumsets}, Graduate Texts in
Mathematics 165, Springer-Verlag, New York, 1996.

\bibitem{Frob} J. L. Ram\'{i}rez Alfons\'{i}n,
\emph{The Diophantine Frobenius problem},
Oxford Lecture Series in Mathematics and its Applications, 30,
Oxford University Press, Oxford, 2005.

\bibitem{Stanchescu-3k-4}Y.~Stanchescu, On addition of two distinct sets of integers,
\emph{Acta Arith.}  \textbf{75}  (1996),  no. 2, 191--194.

\bibitem{taobook} T.~Tao and V.~Vu, \emph{Additive Combinatorics},
Cambridge Studies in Advanced Mathematics 105, Cambridge
University Press, Cambridge, 2006.




\end{thebibliography}
\end{document}